\documentstyle[12pt]{article}
\newtheorem{theorem}{Theorem}\newtheorem{lemma}{Lemma}
\newtheorem{corollary}{Corollary}
\newtheorem{conjecture}{Conjecture}

\newtheorem{proposition}{Proposition}
\newtheorem{example}{Example}

\newtheorem{claim}{Claim}
\newtheorem{problem}{Problem}

\newcommand{\text}[1]{\quad\mbox{#1}\quad}
\def\beq{\begin{equation}}\def\eeq{\end{equation}}
\def\beqn{\begin{eqnarray}}\def\eeqn{\end{eqnarray}}
\def\pont{\hspace{-6pt}{\bf.\ }}

\def\qed{\ifhmode\unskip\nobreak\fi\quad\ifmmode\Box\else$\Box$\fi}

\title{Domination in transitive colorings of tournaments}

\author{
D\"om\"ot\"or P\'alv\"olgyi\footnote{Research supported by Hungarian National Science Fund (OTKA), grant PD 104386
and NN 102029 (EUROGIGA project GraDR 10-EuroGIGA-OP-003), and the J\'anos
Bolyai Research Scholarship of the Hungarian Academy of Sciences.}\\
\small Computer Science Department\\[-0.8ex]
\small Institute of Mathematics\\[-0.8ex]
\small E\"otv\"os Lor\'and University\\[-0.8ex]
\small P\'azm\'any P\'eter s\'et\'any 1/c\\[-0.8ex]
\small Budapest, Hungary, H-1117\\[-0.8ex]
\small \texttt{dom@cs.elte.hu} \and
Andr\'as Gy\'arf\'as\footnote{Supported by Hungarian National Science Fund (OTKA) K104343.}\\
\small Computer and Automation Research Institute\\[-0.8ex]
\small Hungarian Academy of Sciences\\[-0.8ex]
\small Budapest, P.O. Box 63\\[-0.8ex]
\small Budapest, Hungary, H-1518\\[-0.8ex]
\small \texttt{gyarfas@sztaki.hu}
}

\begin{document}
\maketitle

\begin{abstract}
An edge coloring of a tournament $T$ with colors $1,2,\dots,k$ is
called \it $k$-transitive \rm if the digraph $T(i)$ defined by the
edges of color $i$ is transitively oriented for each $1\le i \le
k$. We explore a conjecture of the second author:

\noindent {\em For each positive integer $k$ there exists a
(least) $p(k)$ such that every $k$-transitive tournament has a
dominating set of at most $p(k)$ vertices}.

We show how this conjecture relates to other conjectures and
results. For example, it is a special case of a well-known
conjecture of Erd\H os, Sands, Sauer and Woodrow \cite{SSW} (so
the conjecture is interesting even if false). We show that the
conjecture implies a stronger conjecture, a possible extension of a result
of B\'ar\'any and Lehel on covering point sets by boxes.  The
principle used leads also to an upper bound $O(2^{2^{d-1}}d\log d)$
on the $d$-dimensional box-cover number that is
better than all previous bounds, in a sense close to best
possible. We also improve the best bound known in $3$-dimensions from $3^{14}$ to $64$ and propose possible further improvements
through finding the maximum domination number over parity tournaments.
\end{abstract}

\section{Introduction.}

\subsection{Two conjectures.}
A digraph is called {\em transitive } if its edges are transitively
oriented, i.e. whenever $ab,bc$ are edges $ac$ is also an edge. An
equivalent definition is obtained if vertices represent elements of
a partially ordered set $P$ and $ab$ is an edge if and only if
$a<_Pb$. If orientations of the edges are disregarded, transitive
digraphs define {\em comparability graphs.}

We call an edge coloring of a tournament $T$ with colors
$1,2,\dots,k$  \it transitive \rm if the digraph $T(i)$ defined by
the edges of color $i$ is transitive for each $i, 1\le i \le k$.
If a tournament has a transitive coloring with $k$ colors, we say
that it is {\em $k$-transitive.} Notice that  a $k$-transitive
tournament $T=\cup_{i\in [k]} T(i)$ remains $k$-transitive if for
some index set $I\subseteq [k]$ we replace $T(i)$ by its reverse
for all $i\in I$. This is called a {\em scrambling} of $T$
and is denoted by $T_I$. Clearly, $1$-transitive tournaments are
equivalent to transitive tournaments, i.e. to acyclic tournaments.
In fact, this is true for $2$-transitive tournaments as well.

\begin{proposition}\pont\label{tr21a} Suppose that a tournament $T$
is transitively $2$-colored. Then $T$ is a transitive tournament.
\end{proposition}

\noindent \bf Proof. \rm Observe that $T$ can not contain a cyclic
triangle because two of its edges would be colored by the same color
and that contradicts transitivity of the digraph of that color. Then
$T$ can not have any oriented cycle either, thus $T$ is transitive.
\qed

A vertex set $S$ of a tournament is dominating if for every vertex
$v\notin S$ there exists some $w\in S$ such that $wv$ is an edge
of the tournament. The size of a smallest dominating set in a
tournament $T$ is denoted by $dom(T)$. A well-known fact (an
illustration of the probability method) is that there are
tournaments $T$ with arbitrary large $dom(T)$. Tournaments that
cannot be dominated by $k$ vertices are sometimes called
$k$-paradoxical. The smallest $k$-paradoxical tournaments for
$k=1,2,3$ have $3,7,19$ vertices, respectively. As in \cite{ABKKW},
we associate a hypergraph $H(T)$ to any digraph $T$, with
vertex set $V(T)$ and with hyperedges associated to $v\in V(T)$
as follows: the hyperedge $e(v)$ contains $v$ and all vertices $w$
for which $wv\in E(T)$. In terms of $H(T)$, a dominating set $S$
is a transversal of $H(T)$ and $dom(T)=\tau(H(T))$, the
transversal number of $H(T)$.

In this note we explore two conjectures of the second author and
show that they are equivalent.

\begin{conjecture}\pont\label{transdom} For each positive integer $k$ there is a (least)
 $p(k)$ such that for every $k$-transitive tournament $T$, $dom(T)\le p(k)$.
\end{conjecture}

Conjecture \ref{transdom} is open for $k\ge 3$. In Section
\ref{tropaley} we give some examples of $3$-transitive tournaments
and prove that Paley tournaments cannot provide counterexamples to
Conjecture \ref{transdom}: for every $k$,  Paley tournaments of
order at least $2^{4k+4}$ are not $k$-transitive tournaments
(Theorem \ref{paleynotrans}).

%Atirtam, a regi def szerint 2p(k) is jo korlat lett volna mar
Conjecture \ref{transdom} relates to other conjectures and
results. A set $S\subset V(T)$ is an {\em enclosure set} in
a $k$-colored tournament $T=\cup_{i\in [k]} T(i)$ if for
any $b\in V(T)\setminus S$ there exist an $i\in [k]$ and $a,c\in S$ such that
$ab,bc\in E(T(i))$. We say in this case that $b$ is {\em between $\{a,c\}$}. The
smallest enclosure set of a tournament $T$ is denoted by
$encl(T)$.

\begin{conjecture}\pont\label{encl} For each positive integer $k$ there is a (least)
 $r(k)$ such that every for $k$-transitive tournament $T$, $encl(T)\le r(k)$.
\end{conjecture}

\subsection{Conjecture \ref{transdom} implies Conjecture
\ref{encl}.}

Notice that by definition $p(k)\le r(k)$. Our main observation is
that $r(k)$ is bounded in terms of $p(k)$, i.e. Conjecture
\ref{encl} follows from Conjecture \ref{transdom}.
Denote by $T_I$ the $2^k$ possible scramblings of a
$k$-transitive tournament $T$.

\begin{theorem}\pont\label{sum}$r(k)\le \sum_{I\subseteq [k]} dom(T_I)$.
\end{theorem}

If we know that $p(k)$ is finite, Theorem \ref{sum} implies

\begin{corollary}\pont\label{domencl}$r(k)\le
2^{k}p(k)$.
\end{corollary}

\noindent \bf Proof of Theorem \ref{sum}. \rm
Select a dominating set $S_I$ from each scrambling $T_I$.
Clearly, it is enough to prove the following.

\begin{claim}\pont\label{union} $R=\cup_{I\subseteq [k]} S_I$ is an enclosure set of
$T$. \end{claim}

Indeed, suppose $v\in V(T)\setminus R$ and $v$ is not between any two points of $R$.
This means that for every color $i\in [k]$ either nobody dominates $v$ from $R$,
or nobody is dominated by $v$ from $R$ (or possibly both).
Denote the latter set of colors by $I\subset [k]$.
In this case nobody dominates $v$ in $T_I$, which contradicts the choice of $S_I$.\hfill { } \qed

\subsection{$d$-coordinate tournaments, improved bound on the box-cover number.}

Our main result comes by applying Corollary \ref{domencl} to
special transitively colored tournaments $T$, where $dom(T)$ is bounded. These tournaments are defined by finite sets
$S\subset R^d$. The tournament $T$ has vertex set $S$, the edge
between $p=(p_1,\dots,p_d)$ and $q=(q_1,\dots,q_d)$ is oriented from
$p$ to $q$ if $p_1<q_1$
%$\min\{p_1,q_1\}$ to $\max\{p_1,q_1\}$
(we shall assume that coordinates of points of $S$ are all different).
The color of an edge is assigned according to the $2^{d-1}$ possible
relation of the other $d-1$ coordinates of $p$ and $q$. Thus we have
here $2^{d-1}$ colors and clearly each determines a transitive digraph,
thus $T$ is a $2^{d-1}$-transitive tournament. We call tournaments obtained
from $T$ by scrambling (i.e. reversing direction on any subset of
colors) {\em $d$-coordinate tournaments}. Note that $2^{2^{d-1}}$
$d$-coordinate tournaments are defined by any $S\subset R^d$.

Set $g(d)=\max_{T} encl(T)$ over all $d$-coordinate
tournaments $T$. According to Conjecture \ref{encl}, $encl(T)\le
r(2^{d-1})$. We shall describe this problem in a more geometric
way as follows.

Given two points $p,q\in R^d$, define $box(p,q)$ as the smallest
closed box whose edges are parallel to the coordinate axes and
contains $p$ and $q$. So $box(p,q)=\{x\in R^d\mid \forall i\;
\min(p_i,q_i)\le x_i\le \max(p_i,q_i)\}$.

Denote by $g(d)$ the smallest number such that we can select a set
$P, |P|\le g(d)$ from any finite collection $S\subset R^d$, such
that for any $s\in S$ we have $p,q\in P$ for which $s\in box(p,q)$.
Note that it is possible to give $2^{2^{d-1}}$ points in $R^d$ so
that none of them is in the box generated by two others, implying
$2^{2^{d-1}}\le g(d)$. This lower bound comes from many sources, in
fact related to Ramsey problems as well, see \cite{BL}. In fact, it
is sharp: in any set of $2^{2^{d-1}}+1$ points of $R^d$ there is a
point in the box generated by two of the other points. This nice
proposition comes easily from iterating the Erd\H os - Szekeres
monotone sequence lemma. (A proof is in \cite{P} while \cite{ABKKW}
refers to this as an unpublished result of N. G. de Bruijn.)

On the other hand, B\'ar\'any and Lehel \cite{BL} have shown that
$g(d)$ always exists and gave the upper bound $g(d)\le
(2d^{2^d}+1)^{d2^d}$, which was later improved to $g(d)\le
2^{2^{d+2}}$ by Pach \cite{P}, and finally to $g(d)\le
2^{2^d+d+\log d+\log\log d+O(1)}$ by Alon et al. \cite{ABKKW}.
Note that this last bound is about the square of the lower bound.
We use Corollary \ref{domencl} to improve their method from
\cite{ABKKW} and this gives our main result, an upper bound that
is very close to the lower bound.

\begin{theorem}\pont\label{improve} $g(d)= O(2^{2^{d-1}}d\log d)$. \end{theorem}

We give the proof of this theorem in Section \ref{sec2}.

\subsection{$3$-coordinate tournaments and Parity tournaments.}

Although Theorem \ref{improve} is ``close'' to the lower bound $2^{2^{d-1}}$, it is not useful for small $d$. Note that $g(1)=2$ is obvious, $g(2)=4$ is easy but the best known bound in $3$ dimension was $g(3)\le 3^{14}$ by B\'ar\'any and Lehel.
Here we look at the $3$-dimensional case more carefully and obtain a more reasonable bound of $g(3)$.

In the $3$-dimensional case any point set defines sixteen $3$-coordinate tournaments.
The obtained tournaments are of three kind.

{\bf Case 1.} Six tournaments are ``dictatorships'',
where one coordinate defines the edges and the other two are
irrelevant. These are transitive and thus each is dominated by a
single vertex.

{\bf Case 2.} Eight tournaments are isomorphic to a
$2$-majority tournament, where $uv$ is an edge if two coordinates of
$u$ are bigger than the respective coordinates of $v$. It was proved
in \cite{ABKKW} that such tournaments have a dominating set of three
vertices (and this is sharp).

{\bf Case 3.} The last two tournaments are inverse to each other,
and we call them {\em parity tournaments} as
$uv$ is an edge if an even (resp. odd) number of the three coordinates of
$u$ are bigger than the respective coordinates of $v$. Set $m=\max_{T} dom(T)$ over all parity tournaments $T$.
We strongly believe that $m$ can be determined by a (relatively) simple, combinatorial argument (like in \cite{ABKKW} for $2$-majority tournaments) thus we pose

\begin{problem}\pont Determine $m=\max_{T} dom(T)$ over all parity tournaments $T$.  %Is there a parity tournament which is not a $2$-majority tournament?
\end{problem}

The observations above with Theorem \ref{sum} yield the following upper bound on $g(3)$.

\begin{corollary}\pont \label{3dim}$g(3)\le 6\times 1+8\times 3+2m=30+2m$.
\end{corollary}

It is possible that $m=3$ but we could prove only  $m\le 17$ by a careful modification of the probabilistic proof of Theorem \ref{epsnet} (for details, see the Appendix). Using this and Corollary \ref{3dim} we get $g(3)\le 64$, a first step to a ``down to earth'' bound.

\subsection{A path-domination conjecture and majority
tournaments.}

Conjecture \ref{transdom} also relates to a well-known conjecture of
Erd\H os, Sands, Sauer and Woodrow \cite{SSW}.

\begin{conjecture}\pont\label{ESSW} Is there for each positive integer $k$ a (least)
integer $f(k)$, such that every tournament whose edges are colored
with $k$ colors, contains a set $S$ of $f(k)$ vertices with the
following property:
every vertex not in $S$ can be reached by a monochromatic
directed path with starting point in $S$. In particular, (a.) does
$f(3)$ exist? (b.) $f(3)=3$?
\end{conjecture}

For further developments on Conjecture \ref{ESSW} see \cite{GR,GS,M,MR,PR}.

\begin{proposition}\pont\label{fimpliesp}$p(k)\le f(k)$.
\end{proposition}

\noindent \bf Proof. \rm Suppose $T$ has a transitive $k$-coloring.
From the definition of $f$, there exists $S\subset V(T)$
such that $|S|=f(k)$ and every vertex $w\notin S$
can be reached by a directed path of color $i$ for some $1\le i \le k$
with a starting point $v\in S$. From transitivity of color $i$, $vw\in E(T)$, thus $S$ is
a dominating set in $T$. \hfill { }\qed

Examples \ref{expaley7}, \ref{exblowc3} in Subsection \ref{23} show $p(3)\ge 3$,
thus Proposition \ref{fimpliesp} gives

\begin{corollary}\pont\label{p(3)} If (b) is true in Conjecture \ref{ESSW} then
$p(3)=3$.
\end{corollary}

A last remark relates transitive tournaments to {\it
$k$-majority tournaments $T$}, defined by $2k-1$ linear orders on
$V(T)$ by orienting each pair according to the majority.
Note that every $k$-majority tournament is also a $(2k-1)$-coordinate tournament.
It was
conjectured that $dom(T)$ is bounded by a function  $maj(k)$ for
$k$-majority tournaments and Alon et al. \cite{ABKKW} proved the conjecture with
a very good bound $maj(k)= O(k\log k)$. The existence of $maj(k)$
would also follow from the existence of $p(k)$ (with a very poor bound).

\begin{proposition}\pont$$maj(k)\le p\bigg(\sum_{i=k}^{2k-1} {2k-1\choose
i}\bigg).$$
\end{proposition}

\noindent {\bf Proof.} If $T$ is defined by the majority rule with
linear orders $L_1,\dots, L_{2k-1}$ then one can color each edge
$xy\in T$ by the set of indices $i$ for which $x<_{L_i}y$. It is
easily seen that this coloring is transitive (using at most
$\sum_{i=k}^{2k-1} {2k-1\choose i}$ colors). \hfill { } \qed

%The proof of Theorem \ref{paleynotrans} uses a quasi random property
%of Paley graphs, and gives the idea how to replace it by the
%Regularity lemma to prove Conjecture \ref{transdom} in asymptotic
%sense {Section \ref{asympt}, (Theorem \ref{asthm}).

\section{Proof of Theorem \ref{improve}.}\label{sec2}

We may suppose that $S\subset R^d$ is in general position in the
sense that no two points share a coordinate. Let $T$ denote the
tournament associated to $S$. To bound $dom(T)$, we follow the
same strategy as in Alon et al. \cite{ABKKW}.

As mentioned in the introduction, a dominating set $S$ is a
transversal of $H(T)$, thus $dom(T)=\tau(H(T))$. We need a bound on
$\tau(H(T))$ in terms of $d$. Let $\tau^*(H(T))$ be fractional
transversal number of $H(T)$.

The following claim is from \cite{ABKKW}, we give a different proof
for it.

%ennek nincs vmi koze a lovasz-spencer-vesztergombihoz?
%a linear discrepancyhez? en nem latom, de lehet.
\begin{claim}\pont[Alon et al.] For any tournament $T$, $\tau^*(H(T))< 2$.
\end{claim}
\noindent {\bf Proof.} Let $T$ be a tournament on $n$ vertices, set $H=H(T)$.
We will prove that the fractional matching number,
$\nu^*(H)$, which equals $\tau^*(H)$, is less than $2$. Denote by $A$
the $n\times n$ incidence matrix of $H$, with rows indexed by
vertices, columns with edges. Vectors are column vectors with $n$
coordinates. Denote by $\bf{j}$ the all $1$ vector, by $J$ the all
$1$ $n\times n$ matrix. Since $\nu^*=\max \{\bf{j}^T\bf{x}\mid
A\bf{x}\le \bf{j}, 0 \le x\le j\}$, let us take an $\bf{x}$ such that $A\bf{x}\le
\bf{j}$. This also implies $\bf{x}^TA^T\le \bf{j}^T$. From this,
using that $A+A^T$ contains $1$-s except along the main diagonal,
where all elements are $2$-s, we have
$$({\bf j }^T {\bf
x})^2={\bf{x}}^T{\bf{j}}{\bf{j}}^T{\bf{x}}={\bf{x}}^TJ{\bf{x}}<{\bf{x}}^T(A+A^T){\bf{x}}\le
{\bf{j}}^T{\bf{x}}+{\bf{x}}^T{\bf{j}},$$
implying ${\bf{j}}^T{\bf{x}}<2$. \hfill { }\qed

As in \cite{ABKKW}, we bound the VC-dimension of $H$, $VC(H)$,
then we can use the following consequence of the Haussler-Welz theorem,
formulated first by Koml\'os, Pach and Woeginger, \cite{KPW}.

\begin{theorem}\label{epsnet}\pont $\tau(H)= O(VC(H)\tau^*(H)\log \tau^*(H))$.
\end{theorem}

Note that this implies $dom(T)= O(VC(H(T)))$ in our case. Our next
statement is a slight modification of a respective claim from
\cite{ABKKW}.

\begin{claim}\pont\label{vc} If $T$ is a $d$-coordinate tournament and $h=VC(H(T))$, then $(h+1)^d \geq 2^h$
and so $h\le (1+o(1))d\log d$.
\end{claim}
\noindent {\bf Proof.} Fix $h$ vertices that are shattered. These
vertices divide each of the $d$ coordinates into $h+1$ parts. For
two points that belong to the same $h+1$ parts in each of the $d$
coordinates, the restriction of their dominating hyperedge to the
$h$ vertices is also the same.
(Note that the $h$ vertices themselves do not belong to any of the $(h+1)^d$ possibilities but they are equivalent to a vertex obtained by moving them slightly up or down in all coordinates.)
Since we have $2^h$ different
restrictions, the bound follows. \hfill { }\qed

Putting all together we have
$$\label{dim}
dom(T)=\tau(H(T))= O(VC(H(T))\tau^*(H(T))\log
\tau^*(H(T)))= O(d\log d).
$$

Applying Corollary \ref{domencl} and the above inequality to the
$2^{d-1}$-transitive tournament $T$, we get

$$encl(T)\le 2^{2^{d-1}}dom(T)= O(2^{2^{d-1}}d\log d),$$

which establishes Theorem \ref{improve}.

\section{$2$- and $3$-transitive tournaments, Paley
tournaments.}\label{tropaley}

\subsection{$2$- and $3$-transitive tournaments.}\label{23}

One can say more than stated in Proposition \ref{tr21a} about
$2$-transitive tournaments. For a permutation $\pi=x_1\dots x_n$
of $\{1,2,\dots,n\}$ define the $2$-colored tournament $T(\pi)$ by
orienting each edge $x_ix_j$ from smaller index to larger index and
coloring it with red (resp.\ blue) if $x_i<x_j$ (resp.\ $x_i>x_j$).

\begin{proposition}\pont\label{tr21b} Suppose that a tournament $T$ is
transitively $2$-colored. Then $T=T(\pi)$ with a suitable
permutation $\pi$ of the vertices of $T$.
\end{proposition}

\noindent \bf Proof. \rm  Follows from a result of \cite{PLE}. \hfill { }\qed

Transitive $3$-colorings can be quite complicated. Here we give some
examples.

\begin{example}\pont\label{exbip} Suppose that the vertex set of $T$ is
partitioned into two parts, $A,B$. $T(1)$ is an arbitrary bipartite
graph oriented from $A$ to $B$. $T(2)$ is the bipartite complement,
oriented from $B$ to $A$. $T(3)$ is the union of two disjoint
transitive tournaments, one is on $A$, the other is on $B$. Note
that $T$ can be dominated by two vertices. However, the hypergraph
$H(T)$ can have arbitrarily large $VC$-dimension.
\end{example}

\begin{example}\pont\label{expaley7}(\cite{TH}) The Paley tournament
$PT_q$ is defined for every prime power $q\equiv -1 \pmod{4}$. Its
vertex set is $GF(q)$ and $xy$ is an edge if $y-x$ is a (non-zero)
square in $GF(q)$. Let $PT_7$ be the Paley tournament on $[7]$ with
edges $(i,i+1),(i,i+2),(i,i+4)$. A transitive $3$-coloring of $T$ is
the following:
$$T(1)=\{(1,2),(1,5),(3,4),(3,5),(3,7),(4,5),(6,7)\},$$
$$T(2)=\{(1,3),(2,3),(2,4),(2,6),(4,6),(5,6),(5,7)\},$$
$$T(3)=\{(4,1),(5,2),(6,3),(6,1),(7,1),(7,2),(7,4)\}.$$
\end{example}

{\em Substitution} of a colored tournament $H$ into a vertex $v$ of
a colored tournament $T$ is replacing $v$  by a copy $H^*$ of $H$
and for any $w\in V(H^*),u\in V(T)\setminus V(H^*)$ the color and
the orientation of $wu$ is the same as of $vu$.

\begin{example}\pont\label{exblowc3}Let $T$ be the tournament on nine vertices obtained
from a $3$-colored cyclic triangle by substituting it into each of
its vertices.
\end{example}

 Observe that examples \ref{expaley7} and \ref{exblowc3} both show $p(3)\ge
3$.

\subsection{Paley tournaments.}

By a theorem of Graham and Spencer, Paley tournaments $PT_q$ are
$k$-paradoxical if $q\ge k^24^k$. Thus Paley tournaments are
potential counterexamples to Conjecture \ref{transdom}. This is
eliminated by the next result.

\begin{theorem}\pont\label{paleynotrans} For $q> 2^{4k+4}$,
$PT_q$ has no transitive $k$-coloring.
\end{theorem}
\noindent \bf Proof. \rm
We shall use the following lemma of Alon (Lemma 1.2 in Chapter 9
of \cite{ASP}, where it is stated only for primes but it
comes from \cite{A} where it is stated for prime powers).

\begin{lemma}\pont\label{alon} Suppose $A,B$ are subsets of vertices in the Paley
tournament $PT_q$, for some prime power $q$. Then
$$|e(A,B)-e(B,A)|\le (|A||B|q)^{1/2}.$$
\end{lemma}

Suppose indirectly that $T=PT_q$ has a transitive $k$-coloring with
colors $1,2,\dots,k$. Set $\nu={1\over 2^{2k+2}}$ and define the
type of a vertex $v\in V(T)$ as a pair of subsets $I(v),O(v)$ where
$I(v),O(v)\subseteq [k]$ and $i\in I(v)$ if and only if $d_i^-(v)\ge \nu q$
 and $i\in O(v)$ if and only if $d_i^+(v)\ge \nu q$. Thus
the type of $v\in V(T)$ determines those colors in which the
indegree or the outdegree of $v$ is large.

\begin{claim}\pont\label{io} There is a vertex $w\in V(T)$ whose
type $(I,O)$ satisfies $I\cap O\ne \emptyset$.
\end{claim}
\noindent {\bf Proof of claim. } There are at most $2^{2k}$ types so
there is an $A\subset V(T)$ such that all vertices of $A$ have the
same type and
\begin{equation}\label{largetype}
|A|\ge {q\over 2^{2k}}.
\end{equation}
If for each $a\in A$, $I(a)\cap O(a)=\emptyset$, then $|E(A)|$, the
number of edges in $A$ can be counted as follows:
$${|A|\choose 2}=|E(A)|= \sum_{a\in A,i\in I(a)} d_i^+(a) +\sum_{a\in A,i\notin
I(a)} d_i^-(a)< \nu q|A|$$ where the inequality comes from bounding each term by $\nu q$ (in the first sum we used that 
$i\in I(a)$ implies $i\notin O(a)$). Therefore $|A|<2q\nu+1$, thus, using
(\ref{largetype}) we obtain that
$${q \over 2^{2k}}\le |A|<2q\nu+1$$ implying
$$2^{-(2k+1)}-(2q)^{-1}<\nu.$$  However, by $q>2^{2k+1}$ (being generous here),
$$\nu=2^{-(2k+2)}=2^{-(2k+1)}-2^{-(2k+2)}<
2^{-(2k+1)}-(2q)^{-1}<\nu$$ contradiction. This finishes the proof
of the claim.\hfill { }\qed

Let $w\in V(T)$ be a vertex given by the claim, and set
$A=N^-_i(w),B=N^+_i(w)$ for some $i\in I(w)\cap O(w)$. From transitivity of
color $i$, all edges of $[A,B]$ are oriented from $A$ to $B$ (in
color $i$). Therefore, using Lemma \ref{alon},
$$|A||B|=|E(A,B)|=|E(A,B)|-|E(B,A)|\le
(|A||B|q)^{1/2}$$ leading to $|A||B|\le q$. However, from the
definition of $A,B$, $\nu^2q^2\le |A||B|$. Therefore we get $q\le
{1\over \nu^2}$ contradicting the assumption $q> 2^{4k+4}$. \hfill { }\qed

%\bf Problem 3. \rm Let $T$ be a tournament whose edges are colored
%with three colors so that every cyclic triangle is colored with at
%most two colors. Does there exist a vertex in $T$ from which every
%other vertex can be reached by a monochromatic directed path?

%\bf Remarks. \rm It was proved in \cite{SSW} that $f(2)=1$ (in fact
%in a more general form),  $f(3)\ge 3$, and mentioned that Problems
%1,2 were also asked by Paul Erd\H os.

%Problem 3 has an affirmative answer if the condition (no cyclic
%$3$-colored triangle) on the coloring is restricted further,
%\cite{M}, \cite{MR}, \cite{PR}, \cite{GR}, \cite{GS}.

%\section{ Proof of asymptotic version of Conjecture
%\ref{transdom}}\label{asympt}

%Theorem \ref{paleynotrans} uses a key property of pseudorandom
%tournaments, Lemma \ref{alon}. Since we learned from Szemer\'edi
%that all graphs are random in a certain sense, there is a
%possibility to prove Conjecture \ref{transdom} in asymptotic sense.
%The directed version of the Regularity Lemma is proved by Alon
%\cite{ASH} and colored tournament version is also available,
%\cite{AM}.

%Thus we can start from an $\epsilon$-regular partition of a large
%enough transitively $k$-colored tournament into at least $clno(k)$
%clusters.

\subsubsection*{Acknowledgment}
The authors are grateful to the anonymous referee for the very careful reading and useful suggestions.

\section*{Appendix - Domination in parity tournaments}
Here we prove that $m=\max_{T} dom(T)\le 17$, where the maximum is over all parity tournaments $T$. The proof closely follows the proof of Theorem \ref{epsnet}, for details, see e.g. \cite{Mbook}.

Denote by $H$ the hypergraph associated to some parity tournament.
We have that $\tau^*(H)<2$ and $\pi(n)\le (n+1)^3$ where $\pi(n)$ is the {\em shatter function} of $H$, defined as $$\max_{|S|=n} |\{S\cap F\mid F\in H\}|.$$
First we prove $m\le 19$.

Pick a multiset of $38$ points at random according to the distribution given by $\tau^*(H)$ and divide it randomly into two sets, $A$ and $B$, with $a=19$ and $b=19$ points.
Denote by $E$ the event that there is an $F\in H$ which is disjoint from $A$ and contains $B$.
The probability that $A$ is not a $\frac 12$-net is at most $2^{b}$ times the probability of $E$.
For any given $F$, the probability of $F\cap A=\emptyset$ and $B\subset F$ is at most $\frac{1}{{a+b\choose b}}$.
As $|F\cap (A\cup B)|$ can take at most $\pi(a+b)$ values, the probability that $A$ is not a $\frac{1}{2}$-net is at most $\frac{2^{b}\cdot (a+b+1)^3}{{a+b\choose b}}=\frac{2^{19}\cdot 39^3}{{38\choose 19}}<1$.

Now this bound can be further improved using the fact that instead of working with $\pi$, for us it is enough to look at the subsets $F\cap (A\cup B)$ whose size is exactly $b$.
We define $\pi_k(n)$ as $$\max_{|S|=n} |\{S\cap F\mid |S\cap F|=k, F\in H\}|.$$
The above argument works if $\frac{2^{b}\cdot \pi_b(n)}{{a+b\choose b}}<1$, which is equivalent to $\pi_b(n)<{a+b\choose b}/2^b$. Now we need some good bounds for  $\pi_k(n)$.

A trivial observation is that $\pi_k(n)$ is upper bounded by $$\sum_{i\equiv k \textit{ mod }2} \pi_i(n)\le \lceil (n+1)^3/2\rceil.$$ (This is already sufficient to show $m\le 18$.) We can further improve this in the following way. Consider the three times $n+1$ intervals to which the $n$ points divide the three coordinates. We say that the {\em rank} of a point $p$ is $(x,y,z)$, if there are $x$, $y$ and $z$ points under it in the respective coordinates. So we have $0\le x,y,z\le n$ and our first observation was that $x+y+z\equiv k\textit{ mod }2$ is necessary. It is easy to see that we also need $x+y+z\ge k$ if at least $k$ points are dominated by $p$ and $x+y+z\le 3n-(n-k)$ if at most $k$ points are dominated by $p$. The cardinality of the points not satisfying the first, resp. second inequality, is ${k+2\choose 3}$ and, resp. ${n-k+2 \choose 3}$.

Adding these constraint would make our parity argument a bit more complicated, so instead of the second inequality we only require $x+y+z\le 3n-(n-k)+1$. Using these we get $\pi_{k}\le ((n+1)^3-{k+2\choose 3} -{n-k+1 \choose 3})/2$. If we choose $n=31$ and $k=14$, we get that this is less than ${n\choose k}/2^k$, proving that $m\le n-k=17$.

%http://www.wolframalpha.com/input/?i=n%3D31%2C+k%3D14%2C+k!%28n-k%29!2^k%28%28n%2B1%29^3-%28k%2B2%29%28k%2B1%29k%2F6+-%28n-k%2B1%29%28n-k%29%28n-k-1%29%2F6%29%2F%282n!%29
\end{document}